# WAVELET-BASED ESTIMATION WITH MULTIPLE SAMPLING RATES

By Peter Hall and Spiridon Penev

*Australian National University and University of New South Wales*

We suggest an adaptive sampling rule for obtaining information from noisy signals using wavelet methods. The technique involves increasing the sampling rate when relatively high-frequency terms are incorporated into the wavelet estimator, and decreasing it when, again using thresholded terms as an empirical guide, signal complexity is judged to have decreased. Through sampling in this way the algorithm is able to accurately recover relatively complex signals without increasing the long-run average expense of sampling. It achieves this level of performance by exploiting the opportunities for near-real time sampling that are available if one uses a relatively high primary resolution level when constructing the basic wavelet estimator. In the practical problems that motivate the work, where signal to noise ratio is particularly high and the long-run average sampling rate may be several hundred thousand operations per second, high primary resolution levels are quite feasible.

**1. Introduction.** In this paper we suggest methods for online signal recovery, when a noisy signal is sampled at discrete times and either the raw data, or an estimator of the signal computed from the raw data, is recorded or transmitted after a relatively short time delay. There may be no opportunity to go back and re-sample the signal if it transpires that parts of the signal are so complex that insufficient information was acquired in the first sampling operation. However, there is a possibility of increasing the sampling rate online, if, at the current sampling time, it appears that the rate is insufficient to capture important features of the signal. Nevertheless, the long-run average cost of sampling, per unit time, should not exceed a given bound, imposed (e.g.) by the capacity of the storage device.

How should we shift from one sampling rate to another, and back again? What sorts of gains in performance can we expect to achieve using this technology? In the setting of wavelet estimators, we suggest an answer to the first of these questions; and, for our particular rate-switching algorithm, we answer the second question. Our results help to underpin recent accounts of this type of methodology, discussed by, for example, Chen, Itoh and Shiki (1998), Aldroubi and Gröchenig (2001) and Anon (2002). Our main arguments and results may be summarized in elementary language, even though their detailed description requires somewhat intricate theory. We give the summary below.

A wavelet estimator is, of course, particularly good at recovering complex signals from noisy data. Nevertheless, if the sampling rate is only $\rho$, then a wavelet estimator is unable to adequately approximate a signal whose frequency approaches $\rho$, particularly if that frequency is only exhibited over a short time interval. Moreover, if the sampling rate is only $\rho$, then we may not even be aware that signals with frequency greater than $\rho$ are present. However, a sudden increase in frequency, at a level somewhat below the "base" sampling rate, say $\rho = \rho_1$, might be interpreted as suggesting that higher frequencies are present. Hence, an increase in resolvable frequencies might reasonably be used to trigger an increase in the sampling rate to $\rho_2$, say.

In a wavelet estimator, high-frequency terms are incorporated after an empirical assessment, based on a threshold, of whether or not the coefficients of those terms are significantly different from zero. We use the occurrence of one or more relatively large values among the coefficients to indicate the presence of high-frequency oscillations, and to trigger an increase in sampling rate, from $\rho_1$ to $\rho_2$. Likewise, the absence of large coefficients among these terms is used to trigger a return to rate $\rho_1$.

This algorithm has a number of variants, including (e.g.) using majority-type rules, applied to sets of resolution levels, to define triggers for increasing or decreasing the sampling rate, and using more than two different sampling rates. An explicit restriction on the amount of time during which we sample at the higher rate can be imposed to ensure a relatively early return to rate $\rho_1$ when there is a danger of exceeding data storage capacity.

Of course, if the long-run sampling cost is kept fixed, then increasing the sampling rate in some parts of the signal inevitably involves reducing it in others, relative to the rate that would be employed if sampling were uniform. This necessarily reduces the fidelity of the signal estimate there, measured mathematically in terms of $L_p$ distance, for example. However, in the context of machine-recorded data that motivates our work, error variances are usually small, and so the $L_p$ error penalties incurred by slightly reducing the sampling rate in places where the signal is relatively "uninteresting" are not likely to be high.

The potential gain is that relatively high frequencies that would normally be overlooked can now be recovered. Provided the time periods where these frequencies occur are relatively short in duration, and assuming our algorithm adapts sufficiently quickly, using a higher sampling rate there will require only a modest reduction in the rate at the more common places where the signal is "quiet."

There is, of course, a vast and rapidly growing literature on statistical properties of wavelet methods. We mention here only the papers of Donoho and Johnstone (1994, 1995) and Donoho, Johnstone, Kerkyacharian and Picard (1995); further literature will be discussed in Section 2. Methods for optimal design, in the setting of wavelet methods, have been suggested by Herzberg and Traves (1994) and Oyet and Wiens (2000), although not in the context treated in the present paper.

## 2. Methodology.

2.1. *Model for data generation.* In practice, while digitally recorded data might be the result of sampling at unequally spaced times, they would, nevertheless, involve sampling at times on a grid. (Not all grid points need have data associated with them, however.) Reproducing an approximation to the true signal may involve a mixture of imputation, to estimate the signal at grid points where it was not sampled, and interpolation or smoothing, to reduce the impact of noise. If the edge length of the grid is $\xi$, then the grid points are $k\xi$, for $-\infty < k < \infty$.

Usually, $\xi$ would equal the minimum possible spacing between adjacent sampling times $T_i$, which are indexed in increasing order, are distinct, are integer multiples of $\xi$, and (conceptually) increase from the infinite past to the infinite future. At time $T_i$ we record datum $Y_i$, given by

$$(2.1) \qquad\qquad Y_i = g(T_i) + \varepsilon_i, \qquad -\infty < i < \infty,$$

representing the value of the true signal $g$ at time $T_i$, degraded by additive noise $\varepsilon_i$. The $\varepsilon_i$'s are assumed to have zero mean and variance $\sigma^2$, and the threshold for the wavelet estimator will be constructed so that it is proportional to $\sigma$ and inversely proportional to the square root of sampling rate. In this way it will reflect signal-to-noise ratio.

The value of $T_i$ is determined by previous data, and so is measurable in the sigma-field $\mathcal{F}_{i-1}$ generated by the set of pairs $(T_j, Y_j)$ for $j \leq i-1$. It will be assumed that the distribution of $\varepsilon_i$, conditional on $\mathcal{F}_{i-1}$, has zero mean and finite variance and does not depend on $i$.

2.2. *Data imputation.* At any grid point $t = k\xi$, not necessarily one of the $T_i$'s, we wish to estimate $g(t)$ using the data at the possibly unequally spaced times $T_i$. There is a range of ways of achieving this goal, adapted (for example) from methods suggested for dealing with nonregularly spaced design in more conventional problems where wavelet estimators are employed. See, for example, Hall and Turlach (1997), where interpolation is suggested; Cai and Brown (1998) and Hall, Park and Turlach (1998), where transformation and binning are employed; Cai and Brown (1999), who used a universal thresholding technique; Sardy, Percival, Bruce, Gao and Stuetzle (1999), who considered a variety of different methods; Zhang and Zheng (1999), who addressed theoretical issues associated with nonregular design; Kovac and Silverman (2000), who discussed coefficient-dependent thresholding; Antoniadis and Fan (2001), who developed a penalization approach; Delouille, Franke and von Sachs (2001), Delouille, Simoens and von Sachs (2001) and Delouille and von Sachs (2002), who introduced "design-adapted" wavelet methods for a variety of applications; and Pensky and Vidakovic (2001), who described theory for projection-based techniques.

While these methods have excellent numerical and theoretical properties, and can be expected to produce very good results in relatively familiar settings, not all are suitable for online applications. This is particularly true of relatively computer-intensive techniques, and of those that require that an overview be taken of the full design distribution before determining how the final estimator will be constructed. We shall borrow from Hall and Turlach (1997) and Hall, Park and Turlach (1998), and at each grid point $k\xi$ impute a datum $Z_k$, taking it to equal $Y_i$, where $T_i$ is chosen to be as close as possible to $t$ subject to not exceeding $k\xi$. Define $Z^t(s)$ to equal $Z_k$ for $k\xi \le s < \min\{(k+1)\xi, t\}$ and to equal 0 otherwise. The superscript $t$ indicates that the function $Z^t$ is based only on data that are sampled up to time $t$.

To fully appreciate the role played by imputation, it is important to realize that the wavelet estimator will most likely appear only in the very last step of the chain: "record data–store/transmit data–recover signal." It is at this final stage that imputation will occur, well after any decision has been taken about what to store or transmit. Therefore, although it might appear as though some sort of "internal sampling" at the higher sampling rate might avoid the need to interpolate, that will seldom be possible.

2.3. *Wavelet estimator.* Let $\phi$ and $\psi$ denote the "father" and "mother" wavelet functions, respectively, and let $r \ge 1$ be the least integer such that $\int u^r \psi(u)\, du \ne 0$. That is, $\psi$ is of order $r$. Write $p$ for the primary resolution level, put $p_i = 2^i p$ for $i \ge 0$, and define $\phi_j(t) = p^{1/2}\phi(pt - j)$ and $\psi_{ij}(t) = p_i^{1/2}\psi(p_i t - j)$.

The wavelet coefficients for $g$ are $b_j = \int g\phi_j$ and $b_{ij} = \int g\psi_{ij}$, and the corresponding expansion of $g$ is

$$(2.2) \qquad g = \sum_j b_j\, \phi_j + \sum_{i=0}^{\infty} \sum_j b_{ij}\, \psi_{ij}.$$

The estimators of $b_j$ and $b_{ij}$, based on data observed up to time $t$, are $\hat{b}_j^t = \int Z^t \phi_j$ and $\hat{b}_{ij}^t = \int Z^t \psi_{ij}$, respectively. The hard-thresholded form of our wavelet estimator of $g$ is

$$(2.3) \qquad \hat{g}^t = \sum_j \hat{b}_j^t \phi_j + \sum_{i=0}^{q-1} \sum_j \hat{b}_{ij}^t I(|\hat{b}_{ij}^t| \geq \delta) \psi_{ij},$$

where $\delta > 0$ denotes the threshold, and $q > 1$ needs to be chosen. A soft-thresholded estimator may be constructed similarly.

We shall take

$$(2.4) \qquad \delta = C\sigma(\rho^{-1} \log \rho)^{1/2},$$

where $C > 0$ is a constant. This choice reflects the fact that the variance of $\hat{b}_{ij}^t$ is, in the case of truly independent errors and to a good first approximation, proportional to $\sigma^2/\rho$. Owing to the sequential nature of sampling, the errors are not actually independent, but $\hat{b}_{ij}^t$ is, nevertheless, a martingale, and using that result the same variance approximation can be derived.

In practice, $C$ and $\sigma^2$ usually would be chosen through prior experience with both signals and equipment. In the type of application we envisage, there would be no opportunity for a technician to adjust algorithms *in situ*; the equipment would be expected to function as a "black box." Therefore, the only options are fixed, prior choice of parameters, or automatic, locally adaptive choice. Arguments based on the needs for robustness, real-time analysis and computational economy, and the fact that traditional measures of performance do not apply in this problem, relegate in favor of the latter approach.

It is conventional to take $C = 2$ in the threshold, although $C > 2^{1/2}$ is adequate for our purposes, as we shall show in Section 4. In the case of heavy-tailed data one could use a larger moderate-deviation compensator than the factor $(\log \rho)^{1/2}$ that we employ in (2.4). The compensator $\log \rho$ is sometimes suggested as an alternative.

2.4. *Time delays in near-real time inference.* A feature that distinguishes the context of the present paper from more conventional curve-estimation problems is its "online, real-time" nature. There are two aspects to this, arising when recording and "playing back" the data, respectively. When recording the data we wish to detect sharp increases or decreases in frequency relatively quickly, and to change the sampling rate accordingly. Here, the time-delay should ideally be very small.

The recorded data might, for example, represent acoustic information stored on a CD track that we have played up to time $t$. We wish to produce an approximation to the signal at $t$. When playing the data back in this way it will usually not be a problem if we interpret "at $t$" a little liberally. For instance, it is not essential that the sound we hear at $t$ represent the signal which, at that very time, is being uptaken from the CD by the laser reader. We are prepared to accept a short time delay, the length of which is not as crucial as in the recording phase. In other settings, for example, where the recorded data represent remotely sensed information that will be subjected to detailed analysis in a laboratory, time delay at playback will be even less of an issue.

As we shall show in a moment, it is convenient to take time delays at both recording and playback stages to be inversely proportional to the primary resolution level, $p$, which is an increasing function of the long-term average sampling rate. Hence, the high sampling rates of contemporary digital equipment lead to low time-delays. Moreover, the low noise that often characterizes machine-recorded data, and which permits even larger primary resolution levels, further reduces time-delay.

We shall assume $\phi$ and $\psi$ are continuous on the real line and have compact support, contained within the finite interval $[a, b]$, where $a < 0 < b$. It follows that, for a given index $i$, an integral of the form $\int \alpha \psi_{ij}$, for a function $\alpha$, involves $\alpha(t)$ only if $t \in [(a + j)/p_i, (b + j)/p_i]$. If the integral depends on $\alpha(t)$, then it does not depend on $\alpha(s)$ for $s > t + p_i^{-1}(b - a)$, and, hence, not for $s > t + p^{-1}(b - a)$. A similar argument applies to integrals of the form $\int \alpha \, \phi_i$.

Therefore, if it is acceptable to have a delay of $\tau \equiv (b - a)/p$ time units, between when a signal at $t$ is sampled and when its value at $t$ is estimated, then $t$ can be taken to be an "interior" point of the estimator. In this case estimating $g(t)$ by $\hat{g}^{t+\tau}(t)$ is appropriate. The latter is identical to $\hat{g}^s(t)$ for any $s \geq t + \tau$ and, in particular, is identical to the familiar wavelet estimator that would be used if the full dataset, in infinite time, were employed. It is, therefore, no longer necessary to employ the superscript on $\hat{g}^t$, $\hat{b}_j^t$ and $\hat{b}_{ij}^t$, and we shall usually not use it in the sequel.

The time taken to respond to a change in signal frequency, by increasing or decreasing the sampling rate, will usually equal an integer multiple of $\xi$ that is not less than $\tau$. To appreciate how small $\tau$ might be in practice, note that the optimal choice of $p$, for a high sampling rate, is large. Indeed, the appropriate value of $p^{-1}$ is approximately equal to $(\kappa \sigma^2 / \gamma^2)^{1/(2r+1)} \rho^{-1/(2r+1)}$, where $\sigma^2$ denotes noise variance, $\gamma^2$ is the average of the squared $r$th derivative of the signal, $\kappa$ is a constant depending only on the wavelet type, and $\rho$ (expressed as a frequency in Hz, denoting the number of samples per second) is the sampling rate. See Section 4 for details. Usually, $\kappa^{1/(2r+1)}$ is only a little greater than 1, $\sigma^2$ is small, $\gamma^2$ is moderately large, and $\rho^{-1/(2r+1)}$ is small; see below for discussion. As a result, $\tau$ can be kept to a small fraction of a second.

Sampling rates (and, hence, values of $\rho$) for familiar digital consumer devices are generally quite high. They vary from 8 kHz (for digital telephony), through 32 kHz (digital radio), 44.1 kHz (for conventional CDs) and 96 kHz or 192 kHz for DVD audio, to several mHz for new multi-channel systems. Taking these values to the power $1/5$, so as to model a second-order smoother, gives a small value for $\rho^{-1/5}$. For example, it is 0.1 in the case of a 100 kHz sampler.

2.5. *Rate-switching rule.* Our rule operates on the principle that a signal can be deemed to be relatively erratic, and the sampling rate increased, when the thresholded terms in the double series at (2.3) start to make significant contributions. In theory, the sampling rate can be varied virtually in the continuum. However, we shall treat only a two-rate regime, where the estimator is constructed using rate $\rho_1$ on the majority of occasions, but the rate is increased to $\rho_2$ on relatively rare occasions when high-frequency thresholded terms start to be included in the estimator. Likewise, a reduction in the sampling rate, from $\rho_2$ back to $\rho_1$, is triggered when the threshold inequality $|\hat{b}_{ij}| \geq \delta$ starts to fail to be satisfied.

Therefore, the determination of sampling rate, which is done only at the recording step, uses just the wavelet coefficients $\hat{b}_{ij}$, not the full estimator $\hat{g}$. The latter is employed only at the playback step. Nevertheless, our analysis will treat the two steps together, because the strategy employed during recording must be justifiable by good performance during playback.

2.6. *A specific variable-sampling rate estimator.* If, at the current time, we are sampling the signal at time points that are integer multiples of $\ell\xi$, we shall say we are sampling at rate $\rho_1 = (\ell\xi)^{-1}$. If we are sampling at all integer multiples of $\xi$, we shall say we are sampling at rate $\rho_2 = \xi^{-1}$.

Let $p$ denote a primary resolution level (appearing in the definition of estimators in Section 2.3) that is appropriate when the sampling rate is constant at $\rho_1$ over a long period. Choice of $p$ will be discussed in Section 4. We shall use this $p$ throughout, even when the sampling rate is $\rho_2$, and rely on thresholded terms to produce improved performance when the signal is relatively erratic. However, the values of $q$ at (2.3), and $\delta$ at (2.4), will be rate-dependent. Each will be given the subscript $j$ when the sampling rate is $\rho_j$.

Let $\tau_0$ denote the least integer multiple of $\ell\xi$ that is not less than $(b-a)/p$, where $b-a$ is an upper bound to the widths of the supports of $\phi$ and $\psi$. Then the time-delay $\tau$, introduced in Section 2.4, does not exceed $\tau_0$. Put $q = q_2$ if the sampling rate has been $\rho_2$ for at least the last $\tau_0$ time units, and $q = q_1$ otherwise. Likewise, recalling (2.4) and defining $\delta_j = C\sigma(\rho_j^{-1}\log\rho_j)^{1/2}$, we employ the threshold $\delta_2$ when the sampling rate has been $\rho_2$ for at least the last $\tau_0$ time units, and we use the threshold $\delta_1$ otherwise. In practice, the value of $\sigma$ would be replaced by a value determined after extensive experimentation with real data. Section 4 will discuss choice of $q_1$, $q_2$ and the constant $C$. We could use a smaller time delay than $\tau_0$ when sampling at the higher rate $\rho_2$, but choose not to so as to simplify discussion.

The actual estimator used is given at (2.3). There we take $t = s + \tau_0$ and evaluate the estimator at $s$. It follows that the coefficient estimators $\hat{b}_j$ and $\hat{b}_{ij}$ have the same form they would if the full dataset, in infinite time, were employed. Using this interpretation of $\hat{b}_j$ and $\hat{b}_{ij}$, we define

$$\hat{g}(s) = \sum_j \hat{b}_j \phi_j(s) + \sum_{i=0}^{q-1} \sum_j \hat{b}_{ij} I(|\hat{b}_{ij}| \geq \delta)\psi_{ij}(s),$$

where the rule given in the previous paragraph is used to determine $q$ and $\delta$.

Next we define the mechanism for changing the sampling rate. If we are currently sampling at rate $\rho_1$, then, at time $t = k\ell\xi$, we increase the rate to $\rho_2$ if and only if at least one of the values of $|\hat{b}_{ij}|$, for $j$ such that $\psi_{ij}(t) \neq 0$ and for $i$ such that $p_i$ exceeds a predetermined lower bound $\pi_1$, exceeds the threshold $\delta_1$. If we are currently sampling at rate $\rho_2$, and have been for at least $\tau_0$ time units, we continue at this rate until the next time $t = k\xi$ at which none of the values of $|\hat{b}_{ij}|$, for $\pi_2 \leq p_i \leq p_{q_2}$ and $j$ such that $\psi_{ij}(t) \neq 0$, exceeds $\delta_2$. Here $\pi_2$ is another predetermined lower bound.

Our regularity conditions on $q_1$ and $q_2$ [see (4.3)] do not require $q_1 < q_2$. However, taking $q_1 < q_2$ does reflect the fact that a higher sampling rate allows a greater number of wavelet coefficients to be reliably estimated. Similarly, our assumptions do not demand that $\pi_1 < \pi_2$, but this restriction is not unnatural, for the following reason: $\pi_1$ can be viewed as the highest frequency which the low-sampling-rate estimator is capable of adequately resolving, and $\pi_2$ as the lowest frequency for which sampling at the higher rate is necessary in order to produce an adequate estimate.

Of course, $\delta_1 > \delta_2$, but this does not contradict the fact that exceedences of the thresholds $\delta_1$ and $\delta_2$ are used as parts of rules for increasing and decreasing the sampling rate, respectively. The relatively large size of $\delta_1$ reflects only the fact that sampling at the lower rate produces relatively noisy estimates of wavelet coefficients, which require a relatively high threshold in order to guard against incorrect decisions caused by stochastic variation. It is the values of $\pi_1$ and $\pi_2$, not those of $\delta_1$ and $\delta_2$, which are instrumental in determining whether high- or low-frequency features are present.

Therefore, the rule for switching from rate $\rho_1$ to $\rho_2$ is to increase the rate if and only if

(R.1)  $|\hat{b}_{ij}| > \delta$ for some pair $(i,j)$ with $\psi_{ij}(t) \neq 0$ (at current time $t$) and $\pi_1 \leq p_i \leq p_{q_1}$, where $p = o(\pi_1)$ and $\pi_1 = o(p_{q_1})$;

and the rule for switching back again is

(R.2)  $|\hat{b}_{ij}| \leq \delta$ for each pair $(i,j)$ for which $\psi_{ij}(s) \neq 0$ for some $s \in [t - \tau_0, t]$ (where $t$ denotes current time) and $\pi_2 \leq p_i \leq p_{q_2}$ [where $\rho_1 = o(\pi_2)$ and $\pi_2 \leq p_{q_2}$].

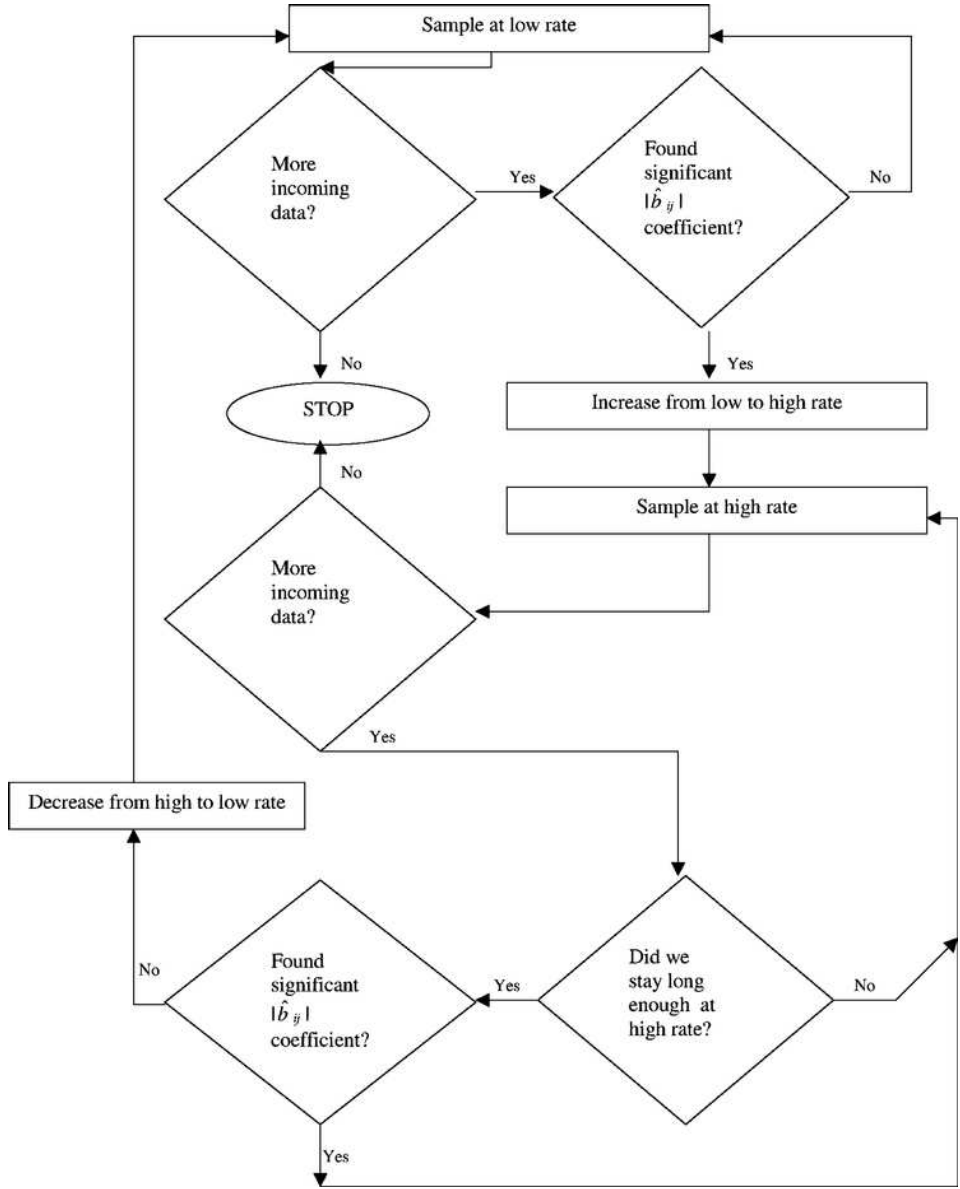

FIG. 1. *Flow chart summarizing rate-switching algorithm.*

[Both (R.1) and (R.2) include regularity conditions which will be used in Section 4.] An overview of the algorithm, after these rate-switching rules are incorporated, is given in the flow chart in Figure 1.

Constraints on the amount of time spent sampling at the higher rate can be introduced to prevent the storage device from filling too rapidly. The algorithm depends on a number of "tuning" parameters, in particular, $p$, $q_1$ and $q_2$, $\pi_1$ and $\pi_2$, the constant $C$ in the threshold formula (2.4), and, of course, the sampling rates themselves. In practice these quantities would be chosen from practical experience with the signal type.

2.7. *Using local Fourier methods with windows of fixed length.* A reviewer has suggested that our wavelet approach might not be competitive relative to a classical local Fourier method using a window of fixed width. However, if, for example, the signal were to have a discontinuity within the interior of the window, or if the values of the signal at either end of the window were unequal, then the Fourier approach—which would perform poorly for functions with discontinuities, interpreted in a periodic sense—would not give good results. In principle this problem could be overcome by choosing the interval adaptively so that discontinuities were situated at its ends. However, that would require continuous local testing for change-points and would arguably be difficult to implement in an on-line fashion. Moreover, such an approach would not address cases where function values at the ends of the interval were different, or where other sorts of signal irregularities, readily adapted to by wavelets, were present.

This issue, of the noncompetitiveness of fixed-bandwidth, local-Fourier methods relative to wavelet ones, in the context of signals with discontinuities and other types of irregularity, is unrelated to our rate-switching scheme. It arises equally in conventional function estimation problems.

## 3. Numerical properties.

3.1. *Smooth signals with aberrations.* Here we illustrate performance in the case of a smooth sinusoid with four different aberration sequences, depicted in Figures 2–5, respectively. Figures 2 and 3 deal with aberrations of increasing amplitude and fixed frequency, and increasing frequency and fixed amplitude, respectively, added at extrema (i.e., at peaks and troughs) of the sinusoid. Figures 4 and 5 address the same respective types of aberration, but added at relatively "linear" places between peaks and troughs. Formulae for the true functions, $g = g_1, \ldots, g_4$, used to produce the respective figures, are given at (3.1)–(3.4).

Figures 3–4 have three panels, showing, respectively, the true signal, its wavelet estimate based on dual-rate sampling and its estimate in the case of fixed-rate sampling using the average of the sampling rates employed in dual-rate sampling. In particular, for a given realization we calculated the number of sampling operations used by the dual-rate algorithm to produce the estimate in the second-to-last panel; and we then sampled at a constant rate, using this number of sampling operations, and employed the data so obtained to produce the estimate in the last panel. Similar results were obtained if, in the constant-rate case, we sampled at the rate obtained by averaging over all $B = 500$ Monte Carlo simulations in the dual-rate case.

The last three panels of Figure 2 show, respectively, the results described above. The first panel of Figure 2 depicts the noisy dataset from which the estimates in the third and fourth panels were computed. We have not shown the noisy data for the other three signals, since doing so adds little of interest.

The superimposed dashed line, in the second-to-last panel of each figure, indicates sampling rate as a function of time. Where the line is at level 0 or 0.5 the sampling rate was low or high, respectively. It can be seen that the rate actually switches up and down several times during high-frequency oscillations. (Using a slightly modified rate-switching rule virtually eliminates these fluctuations and improves performance, but we do not show those results here.)

For each panel of each figure the results shown are those for the realization that gave, among all $B = 500$ realizations of data corresponding to that

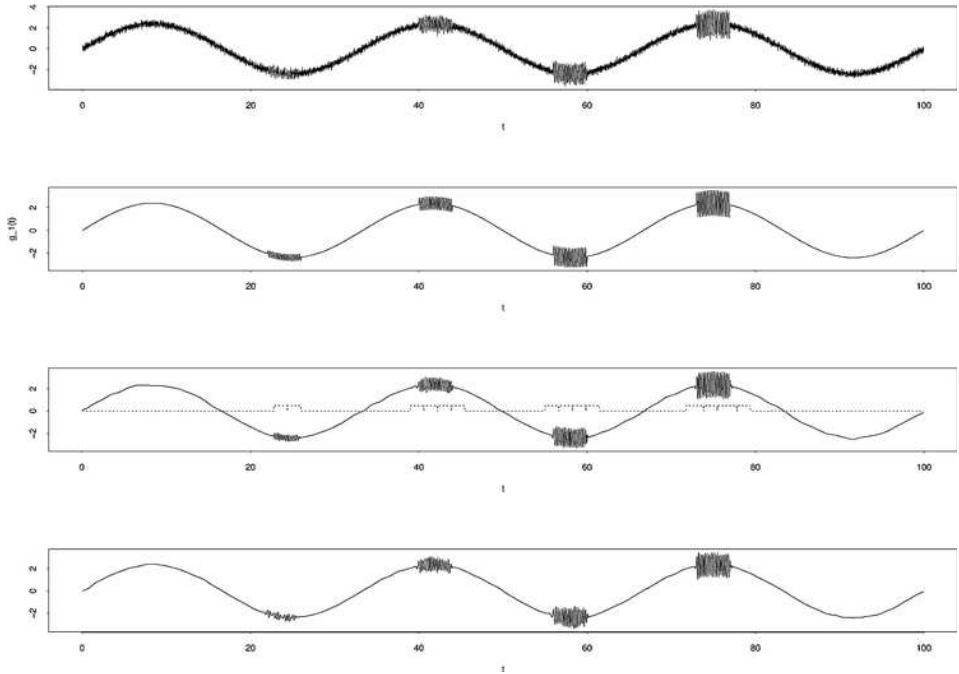

Fig. 2. *Analysis of noisy observations of $g_1$. The first panel shows the noisy data at a sampling rate of 100 Hz, the second panel shows the true signal [i.e., a graph of $y = g_1(t)$, where $g_1$ is given by (3.1)], and the third and fourth panels show the estimates of $g_1$ obtained by dual- and constant-rate sampling, respectively. Here and in subsequent figures, the superimposed dashed line indicates whether the algorithm was operating at the "high" or the "low" sampling rate.*

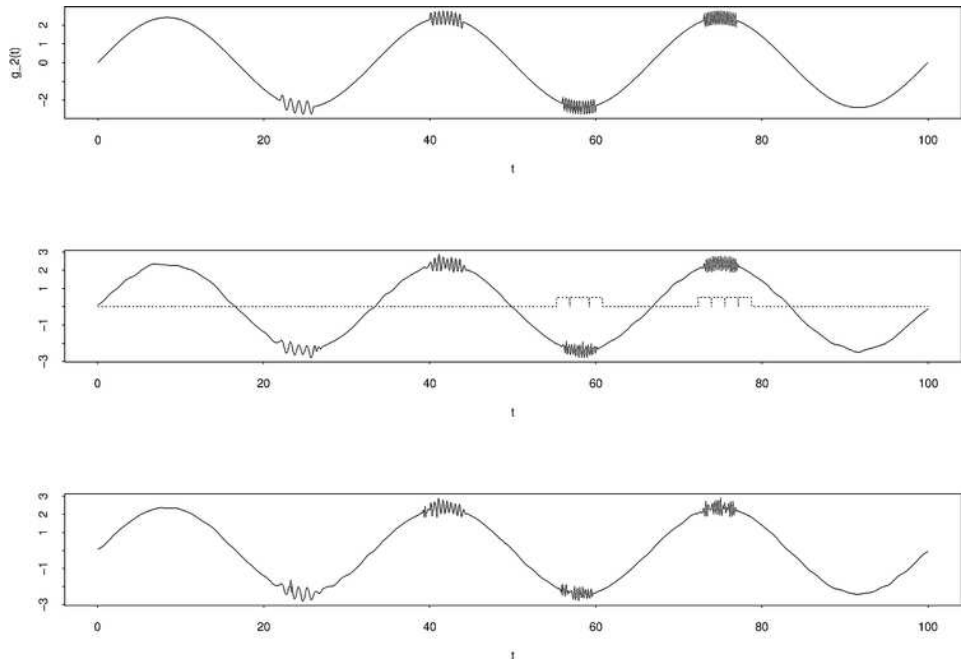

FIG. 3. *Analysis of noisy observations of $g_2$. The first panel shows the true signal [i.e., a graph of $y = g_2(t)$, where $g_2$ is given by (3.2)], and the second and third panels show the estimates of $g_2$ obtained by dual- and constant-rate sampling, respectively.*

figure, the median value of integrated squared error. (For two of the signals we actually conducted $B = 1000$ simulations to check whether the results were significantly different, but they were, in fact, virtually identical. The results reported here are all for the $B = 500$ case.)

Following standard practice we illustrate results in the cases $p = 1$ and $C = 2$, although more favorable results were obtained for different values. In the algorithm discussed in Section 2 we took $q_1 = 4$, $q_2 = 5$, $\pi_1 = 2$ and $\pi_2 = 3$.

The function $\psi$ was chosen from the Daubechies family of compactly supported wavelets with extremal phase and $r = 5$ (i.e., with the length of its support equal to $2r - 1 = 9$).

Signals were sampled at discrete points in the interval $[0, 100]$. To each sampled value, Normal $N(0, 0.15^2)$ noise was added. The edge length of the basic sampling grid (i.e., the minimum permitted spacing between adjacent sampling times) was chosen to be $\xi = 0.01$, which would correspond to a sampling rate of 100 data per unit time if sampling were performed at each grid point. This rate, which we shall refer to as "100 Hz," is the rate used to construct the picture of noisy data in the first panel of Figure 2.

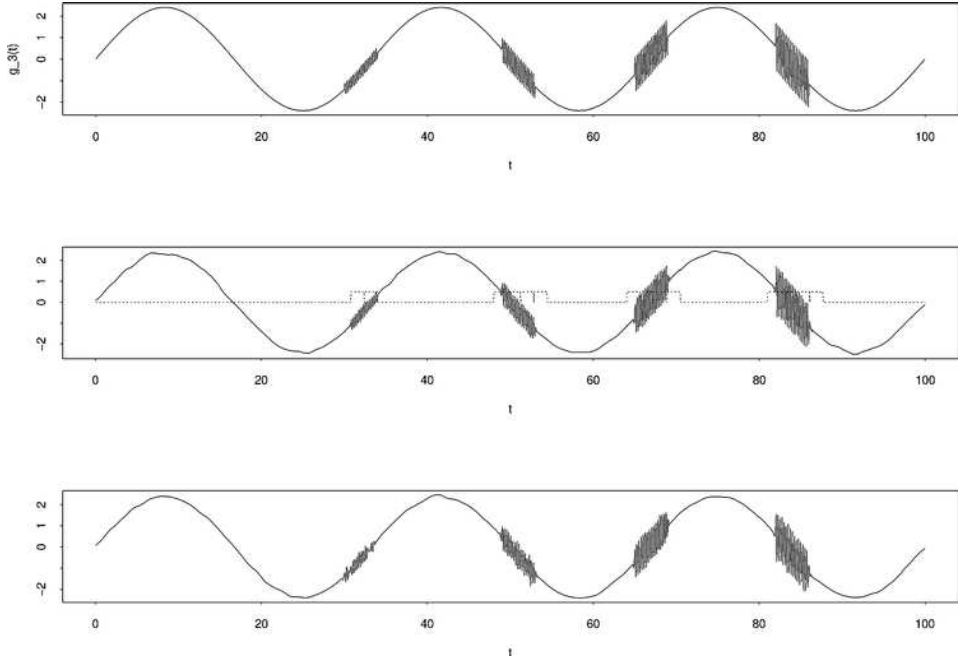

Fig. 4. *Analysis of noisy observations of $g_3$. Panels are in the same order as in Figure* 3.

For our dual-rate algorithm, in the low-rate mode we sampled every sixth observation. That is, in the notation of Section 2.6, we took $\ell = 6$. Equivalently, $\rho_1 = 100/6 \approx 17$ Hz. In the high-rate mode we sampled at each grid point, so that $\rho_2 = 100$ Hz. Sampling at the high rate continued until the end of a minimal time period, of length $(2r-1)100(6p)^{-1} + 1 \approx 150$ units, in which no wavelet coefficient exceeded the threshold.

For each signal type, observations from the first 5% of the time interval $[0, 100]$, sampling at the higher rate, were used to estimate error variance and, thereby, compute thresholds. In particular, we did not assume error variance to be known, although in practice it would most likely be fixed in advance.

Given an interval $[a, b]$, let $I_{[a,b]}(t) = 1$ or 0 according as $t \in [a, b]$ or $t \notin [a, b]$. In this notation, formulae for the functions shown in the first panels of Figures 3–5, and second panel of Figure 2, are, respectively,

$$\begin{aligned}
g_1(t) = {} & 2.4 \sin(0.06\pi t) + 0.2525 \sin\{8\pi(t - 24)\} I_{[22,26]}(t) \\
& + 0.5050 \sin\{8\pi(t - 42)\} I_{[40,44]}(t) \\
& + 0.7575 \sin\{8\pi(t - 58)\} I_{[56,60]}(t) \\
& + 1.1 \sin\{8\pi(t - 75)\} I_{[73,77]}(t),
\end{aligned}$$

(3.1)

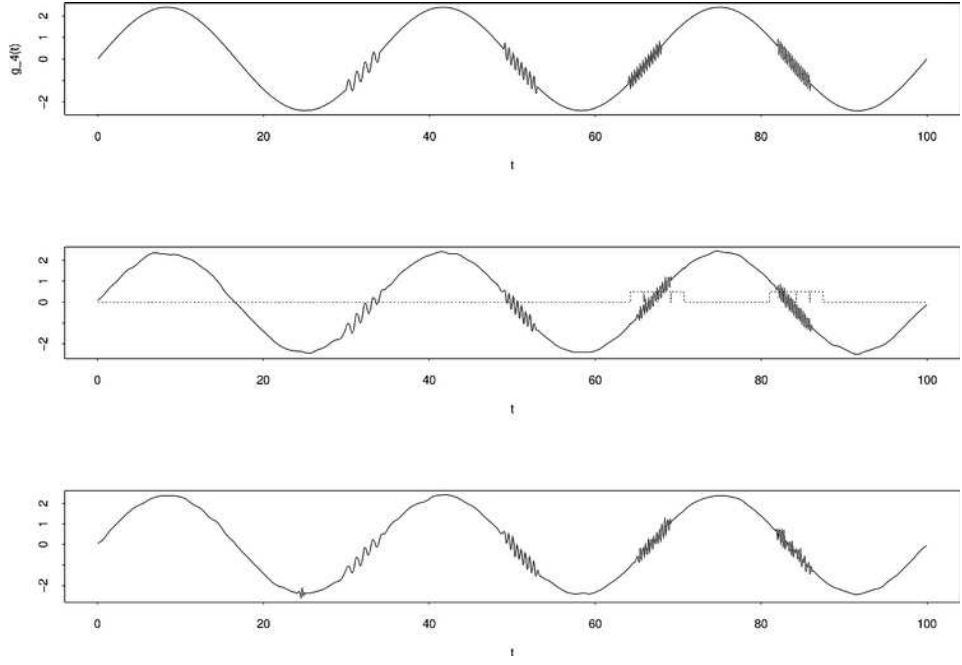

Fig. 5.   *Analysis of noisy observations of $g_4$. Panels are in the same order as in Figure 3.*

$$g_2(t) = 2.4\sin(0.06\pi t) + 0.35\sin\{2\pi(t-24)\}I_{[22,26]}(t)$$
$$+ 0.35\sin\{4\pi(t-42)\}I_{[40,44]}(t)$$
(3.2)
$$+ 0.35\sin\{6\pi(t-58)\}I_{[56,60]}(t)$$
$$+ 0.35\sin\{8\pi(t-75)\}I_{[73,77]}(t),$$

$$g_3(t) = 2.4\sin(0.06\pi t) + 0.2525\sin\{8\pi(t-32)\}I_{[30,34]}(t)$$
$$+ 0.5050\sin\{8\pi(t-51)\}I_{[49,53]}(t)$$
(3.3)
$$+ 0.7575\sin\{8\pi(t-67)\}I_{[65,69]}$$
$$+ 1.1\sin\{8\pi(t-84)\}I_{[82,86]}(t),$$

$$g_4(t) = 2.4\sin(0.06\pi t) + 0.35\sigma\{2\pi(t-32)\}I_{[30,34]}(t)$$
$$+ 0.35\sin\{4\pi(t-51)\}I_{[49,53]}(t)$$
(3.4)
$$+ 0.35\sin\{6\pi(t-67)\}I_{[65,69]}(t)$$
$$+ 0.35\sin\{8\pi(t-84)\}I_{[82,86]}(t).$$

These represent a basic sinusoid, with frequency 0.06 and formula $g(t) = 2.4 \times \sin(0.06\pi t)$, to which are added, in the cases of functions $g_1, \ldots, g_4$, respectively, (i) four aberrations, each with frequency $8\pi$ and with increasing amplitudes, and each of duration four time units; (ii) four aberrations, each with amplitude 0.35 and with increasing frequencies [culminating in the frequency arising in case (i)], and each of duration four time units; and (iii) and (iv) the respective versions of (i) and (ii) where the four aberrations are added midway between extrema of the sinusoid. (The aberrations are of the same duration in each case, although it may appear that durations are longer in the cases of signals $g_3$ and $g_4$.)

The following qualitative properties of dual-rate sampling are illustrated by Figures 2–5. Performance advantages are generally most clear in particularly difficult cases, where the aberrations are of relatively high frequency and low amplitude and so are difficult to distinguish from noise. (The first of the four aberrations added to the signal in Figure 2 is of just this type.) Even though the advantages of dual-rate sampling become more evident as frequency increases, it can, nevertheless, perform well even for relatively low-frequency aberrations (see Figures 3 and 5). Its potential is most marked when an aberration is added to a part of the signal which is changing relatively fast, such as to an extremum of the sine waves (see Figures 2 and 3). However, in difficult cases, where the aberration is of low amplitude and high frequency, it has much to offer in other cases too (see Figures 4 and 5).

The mean integrated squared errors (MISEs) of the four signals, approximated by averaging integrated squared errors over the $B = 500$ simulations conducted for each of the four signals, were (i) 0.348, (ii) 0.320, (iii) 0.341, (iv) 0.321 in the respective cases of $g_1, \ldots, g_4$, for dual-rate sampling; and, respectively, (i) 0.963, (ii) 0.342, (iii) 0.965, (iv) 0.348 for constant-rate sampling with the same average sampling rate. Noting that the MISE advantages of dual-rate sampling are substantially greater in cases (i) and (iii) than in cases (ii) and (iv), one reaches the expected conclusion that dual-rate sampling primarily overcomes problems due to aberrations of a high-frequency, rather than low-amplitude, nature.

Indeed, on the basis of these results one might argue that, in MISE terms, the advantages of dual-rate sampling are marginal in the cases of signals $g_2$ and $g_4$. At first sight this seems at variance with a visual inspection of the figures. However, calculating mean integrated squared errors over only the four intervals, each of length four time units, for each signal, one obtains instead the values (i) 0.202, (ii) 0.180, (iii) 0.185, (iv) 0.162 in the dual-rate case, and (i) 0.869, (ii) 0.245, (iii) 0.763, (iv) 0.209 in the constant-rate setting. Therefore, in the cases of $g_2$ and $g_4$, dual-rate sampling does confer an advantage in terms of its ability to resolve the aberrations, although not as much of an advantage as in the cases of $g_1$ and $g_3$.

3.2. *Discontinuous signals with aberrations.* In order to show that our algorithm is not adversely affected by jump discontinuities in signals, we applied it to signals which had high-frequency aberrations just before, or just after, or shortly before or after, jumps.

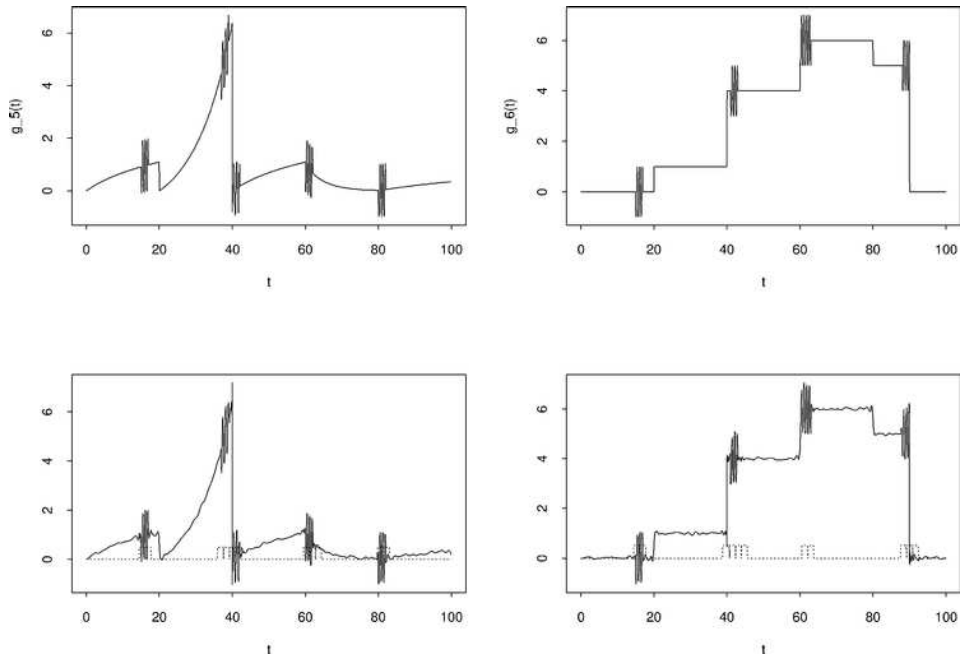

Fig. 6. *Analysis of noisy observations of $g_5$ and $g_6$. The top two panels show the true functions. Immediately below them are the respective estimates, obtained using the rate-switching rule.*

Specifically, the function $g_5$ has jumps at points of increase or decrease in the function, and one jump (the third) at a point which is both a point of increase and a point of decrease. On the other hand, the "block" function $g_6$ has zero derivative except at points that are either part of high-frequency aberrations, or are located at jumps. Apart from the fact that we use different functions in the present section, all settings (and, in particular, all tuning parameters) are the same as in Section 3.1.

Results are summarized in Figure 6, which shows (for each of the two signals) the realization that gave the median value of integrated squared error out of the 500 simulations conducted. In the case of the signal $g_5$, it can be seen that the isolated discontinuity near $t = 20$ causes no problems for the rate-switching rule, and that the method takes in its stride even the very large discontinuity near $t = 40$, which has high-frequency aberrations on either side. The Gibbs phenomena which are present, and which are related to the jump discontinuities, also appeared when our experiments were conducted in the absence of the high-frequency aberrations (but with the jumps still present). These phenomena are a feature of the wavelet method rather than of the rate-switching rule.

The method also gives good results for the function $g_6$, although better performance for such a signal (with or without the high-frequency episodes) usually is obtained using a Daubechies wavelet with a smaller value of $r$. (We employed the same parameter values throughout our study, since we did not wish parameter values to be confounded with function type when interpreting the results.) The isolated discontinuities at $t = 20$ and $t = 80$ are dealt with very well, and the method also enjoys good performance around the points $t = 60$ and $t = 63$. The Gibbs effect to the right of $t = 43$ and $t = 90$ are caused by the large discontinuities there.

Isolated jump discontinuities, such as those in the functions $g_5$ and $g_6$, can sometimes trigger an increase in sampling rate; the jumps may be "misinterpreted" by the algorithm as very high-frequency phenomena. However, this does not cause difficulty. When a sampling-rate increase occurs at a jump discontinuity, it elicits extra information about the location and size of the jump, and that does no harm. Moreover, the rate quickly switches down again after the jump, so little cost is incurred through additional sampling.

For example, in the case of the function $g_6$, jumps at the points 20, 60 and 80 were sufficiently small not to trigger any rate increase. A rate change did generally occur in connection with the larger jump at 40, but without detrimental effects on the estimator. Indeed, the algorithm deduced that the jump was followed by high-frequency events, and correctly maintained sampling at the higher level until the high-frequency sinusoids were past. If the high-frequency events immediately to the right of 40 were removed, then the algorithm returned quickly to the lower sampling rate immediately after the jump.

The functions $g_5$ and $g_6$ are given by the following:

$$\begin{aligned}
g_5(t) = {} & \log(1 + 0.1t)I_{[0,20]}(t) + [\exp\{0.1(t-20)\} - 1]I_{(20,40]}(t) \\
& + \log\{1 + 0.1(t-40)\}I_{(40,60]}(t) + \exp\{-0.2(t-60)\}I_{(60,80]}(t) \\
& + 0.5\log\{1 + 0.05(t-80)\}I_{(80,100)}(t) + \sin\{9(t-16)\}I_{(15,17]}(t) \\
& + \sin\{8(t-38)\}I_{[37,39]}(t) + \sin\{7.3(t-41)\}I_{(40,42]}(t) \\
& + \sin\{9(t-61)\}I_{(60,62]}(t) + \sin\{9(t-81)\}I_{(80,82]}(t),
\end{aligned}$$

$$\begin{aligned}
g_6(t) = {} & I_{(20,40]}(t) + 4I_{(40,60]}(t) + 6I_{(60,80]}(t) + 5I_{(80,90]}(t) \\
& + \sin\{9(t-16)\}I_{(15,17]}(t) + \sin\{9(t-42)\}I_{(41,43]}(t) \\
& + \sin\{9(t-61)\}I_{(60,63]}(t) + \sin\{9(t-89)\}I_{(88,90]}(t).
\end{aligned}$$

**4. Theoretical properties.** The main aim of this section is to establish, under explicit conditions, the four properties below. Together they describe the manner in which the rate-switching algorithm responds to different signal frequencies, and the way in which it can increase the estimator's overall performance. Proofs are given in a longer version of the paper available on the web [Hall and Penev (2002)].

PROPERTY I (Sampling rate remains at $\rho_1$ during "quiet" periods). Suppose we start at the left-hand end of a finite interval $\mathcal{I}$ using sampling rate $\rho_1$, and that the signal is relatively quiet in $\mathcal{I}$. Then the probability that rate $\rho_1$ persists right across $\mathcal{I}$ converges to 1 as $\rho_1 \to \infty$. Moreover, the rate will quickly switch from $\rho_2$ to $\rho_1$ if we start a quiet interval at the higher rate. See Theorems 4.1 and 4.6 for details.

PROPERTY II (Sampling rate increases to $\rho_2$ when signal complexity increases). Suppose the variable-rate estimator is operating at rate $\rho_1$ during a quiet period, and then enters a period of relative activity. Then the algorithm will, with high probability, trigger a switch from rate $\rho_1$ to $\rho_2$ during $\mathcal{J}$. See Theorem 4.4.

PROPERTY III (Sampling rate remains at $\rho_2$ through periods of high-frequency fluctuations). Once the sampling rate has increased from $\rho_1$ to $\rho_2$, it stays there with high probability, provided the signal is sufficiently noisy. See Theorem 4.5.

PROPERTY IV (Dual sampling rates can enhance recovery of high-frequency oscillations, with little adverse affect on estimation of low-frequency features). If sampling is undertaken at rate $\rho_1$, then the estimator is able to recover (in the sense of consistent estimation) signals that have $r$ continuous derivatives, and, indeed, can recover fluctuations that have frequencies of smaller order than $\rho_1$. If sampling is at rate $\rho_2$, then frequencies of smaller order than $\rho_2$ can be recovered. The dual-rate estimator is able to consistently estimate high-frequency parts of the signal that would not be accessible using a constant-rate estimator with the same long-run average sampling cost. This can be done without degrading, in first-order asymptotic terms, the accuracy of approximation in time intervals where the signal is of relatively low frequency. See Theorem 4.7 and the discussion at the end of this section.

Our asymptotic arguments, based on high sampling rates, are justified by the high rates and low noise levels which are commonly encountered in practice; see Section 2.4 for discussion. We shall state our main results, and particularly the regularity conditions, in such a way that the proofs do not require analysis at the level of a martingale error process. Nevertheless, the results have direct application in the latter context, as we shall relate.

Theorems 4.1–4.3 address performance of the wavelet estimator when sampling is at a constant rate and the true signal is smooth. Our aim is to indicate the appropriate sizes for $p$ and $\delta$ in this setting, thereby motivating choices in the variable-sampling rate case when $g$ is not so smooth. Therefore, for the present we take $g = g_0$, where

(4.1)         $g_0$ is an $r$-times continuously differentiable function defined on the whole real line.

Suppose too that

(4.2)         the errors $\varepsilon_i$, in (2.1), are identically distributed with $E|\varepsilon_i|^{6+B} < \infty$ for some $B > 0$, zero mean and variance $\sigma^2$.

Let $\rho$ denote either $\rho_1$ or $\rho_2$. For reasons that will become clear in Theorem 4.2, the appropriate size of $p$ for smooth signals is $\rho^{1/(2r+1)}$, for large $\rho$. To determine the correct size of $q$, observe that we cannot resolve frequencies as large as $\rho$ if we are only sampling at rate $\rho$. Therefore, we shall select $q$ so that $p_q$ is a little smaller than $\rho$; let it be of order $\rho^{1-c}$, where $c > 0$. This is equivalent to $2^q = O(p^{-1}\rho^{1-c})$. Note too that we use the same $p$ when $\rho = \rho_1$ or $\rho = \rho_2$. These considerations motivate the regularity condition

(4.3)         $p = p(\rho_1) \asymp \rho_1^{1/(2r+1)}$ and $q = q(\rho) \to \infty$ so slowly that $2^q = O(p^{-1}\rho^{1-c})$ for some $c \in (0,1)$.

(If $a$ and $b$ are positive functions of $\rho$, then the property $a \asymp b$, as $\rho \to \infty$, means that the ratio $a/b$ is bounded away from zero and infinity along the sequence.)

Finally, suppose that

(4.4)         $\phi$ and $\psi$ are each bounded and supported on the compact interval $[a,b]$, $\psi$ is of order $r$ as defined in Section 2.3, $\int \phi = 1$, and integer translates of $\phi$ are orthonormal.

Our next theorem shows that for functions such as $g_0$, the thresholded terms only very rarely make a contribution to the estimator $\hat{g}$. Recall that $\hat{g}(s) = \hat{g}^{s+\tau_0}(s)$, where $\hat{g}^t$ is given by (2.3). Let $\mathcal{I}$ denote a finite interval.

THEOREM 4.1. *Suppose data are generated by the model at* (2.1), *with signal* $g = g_0$ *and independent errors* $\varepsilon_i$. *Assume the sampling rate is constant at* $\rho = \rho_1$ *or* $\rho_2$, *and that* $\delta$ *in the definition at* (2.3) *is given by* (2.4), *with* $C > 2^{1/2}$. *Suppose too that* (4.1)–(4.4) *hold, with* $B > C^2/(1-c)$ *in* (4.2) *and, in* (4.3), $(q,\rho) = (q_j,\rho_j)$ *in the two respective cases. Then* (a) *if the sampling rate is* $\rho_1$, *the probability that a thresholded term enters nondegenerately into the estimator* $\hat{g}(t)$ *for some* $t \in \mathcal{I}$ *converges to zero as* $\rho_1 \to \infty$;

and (b) *if the sampling rate is $\rho_2 \geq \rho_1$, the probability that a thresholded term corresponding to a resolution level $p_i$ greater than $C_1 \rho_2^{1/(2r+1)}$ (for an arbitrary fixed $C_1 > 0$) enters nondegenerately into $\hat{g}(t)$ for some $t \in \mathcal{I}$ converges to zero as $\rho_2 \to \infty$.*

Our proof of the theorem shows that, for $\rho = \rho_1$ or $\rho_2$, the respective probabilities equal $1 - O(\rho^{1-(C^2/2)+\eta})$ for each $\eta > 0$. This type of bound also applies to all the probabilities that are discussed in Theorems 4.4–4.6: each of the probabilities converges to 1 at rate $\rho^{A(C)}$, where $\rho$ denotes the relevant sampling rate and $A(C)$ can be made arbitrarily large by choosing the constant $C$, in (2.4), sufficiently large.

Part (a) of Theorem 4.1 motivates us to consider in more detail the estimator obtained by sampling at the base rate $\rho_1$. A result in this case was given by [Hall and Patil (1995)](); the following version is better adapted to the present context.

THEOREM 4.2. *Assume data are generated by the model at (2.1), with independent errors $\varepsilon_i$ and time points $T_i$ equally spaced $\rho_1$ units apart. Suppose too that, for $(q, \rho) = (q_1, \rho_1)$, (4.1)–(4.4) hold, and that $\delta$ in the definition at (2.3) is given by (2.4) with $C > 2^{1/2}$. Then, for all finite intervals $\mathcal{I}$,*

$$
\begin{aligned}
(4.5) \qquad & \int_{\mathcal{I}} E(\hat{g} - g_0)^2 \\
& = \rho_1^{-1} p |\mathcal{I}| + p^{-2r} \kappa^2 (1 - 2^{-2r}) \int_{\mathcal{I}} (g^{(r)})^2 + o(\rho_1^{-1} p + p^{-2r})
\end{aligned}
$$

*as $\rho_1 \to \infty$.*

It is immediately clear from (4.5) that for a constant sampling rate $\rho_1$, and a smooth signal $g_0$ that is not a polynomial of degree $r - 1$, the asymptotically optimal value of $p$ will be a constant multiple of $\rho_1^{1/(2r+1)}$. This motivates the first part of (4.3), and suggests that we should take $p$ to be of this size in the variable-sampling rate case too, provided the signal is smooth "most" of the time.

The $L_2$ convergence rate, when $\rho = \rho_1$ and $p \asymp \rho_1^{1/(2r+1)}$, is, as implied by Theorem 4.2, $O(\rho_1^{-r/(2r+1)})$. Based on experience for more conventional estimators, we expect the $L_\infty$ convergence rate to differ from this by no more than a factor $(\log \rho_1)^{1/2}$. Theorem 4.3 confirms this. The reason for our interest in $L_\infty$ rates is that, for more complex signals, we shall use consistency in the supremum metric to assess performance of the estimator.

THEOREM 4.3. *Assume the conditions of Theorem 4.2, and, in particular, that the sampling rates are constant at $\rho_1$. Then,*

$$
\sup_{t \in \mathcal{I}} |\hat{g}(t) - g(t)| = O_p\{\rho_1^{-r/(2r+1)} (\log \rho_1)^{1/2}\}.
$$

Next we develop a model for high-frequency fluctuations. It will be asymptotic in character, and depend on a parameter, $\nu$ say, which we could interpret either as one of the rates $\rho_1$ and $\rho_2$, but perhaps more realistically as the long-term average sampling rate; see (4.12) for a definition of the latter. Our theory will involve $\nu$ diverging to infinity.

Our model for the signal will amount to a smooth function $g_0$, described at (4.1), to which we shall add (on an interval $\mathcal{J}$) fluctuations at least one of which is of unboundedly large frequency. If the frequencies of the fluctuations are represented by $\alpha_1, \alpha_2, \ldots$, then, in order for at least one of them to lead to a rate change as suggested by rule (R.1) in Section 2, we should assume that

(4.6) $\qquad$ for some $k = k(\nu)$, $\qquad \alpha_k / \pi_1 \to \infty \quad$ and $\quad \alpha_k = o(p_{q_1})$.

The high-frequency fluctuations that we shall add to $g_0$ will have the form $\gamma\{\alpha(\cdot - u)\}$, where $\alpha = \alpha_k$ and

(4.7) $\qquad$ $\gamma$ is a nondegenerate function, supported on the interval $[-1, 1]$ and having $r$ continuous derivatives on the real line.

Without loss of generality, $\gamma$ is centred so that

(4.8) $$\gamma^{(r)}(0) \neq 0.$$

(Any shift in the location of $\gamma$ can be incorporated into the $u$'s.) For the sake of simplicity we shall choose $\gamma$ to be the same for each fluctuation, although our results are not changed if we use a more elaborate construction. The locations $u$ and frequencies $\alpha$ will vary, however, as follows.

Since the function $g_0$ satisfies (4.1), then its first $r$ derivatives are bounded in any compact interval. On the other hand, if $\gamma$ satisfies (4.7) and $\alpha = \alpha(n) \to \infty$, then the supremum of the absolute value of any one of the first $r$ derivatives of $\gamma\{\alpha(\cdot - u)\}$ diverges to infinity in any open interval containing $u$. Therefore, $\gamma\{\alpha(\cdot - u)\}$ can fairly be said to exhibit fluctuations whose size is an order of magnitude greater than in the case of $g_0$. We shall use the former function to model high-frequency wiggles which trigger an increase in the sampling rate, from $\rho_1$ to $\rho_2$.

We shall add the fluctuations within an interval $\mathcal{J}$, the length of which could converge to zero as $\nu \to \infty$. Thus, there will be a "cluster of wiggles" $\gamma_k$ in $\mathcal{J}$, described through a sequence of pairs $(u_k, \alpha_k)$ with the following property:

(4.9) $\qquad$ the functions $\gamma_k = \gamma\{\alpha_k(\cdot - u_k)\}$ are all supported in $\mathcal{J}$, and no two of the support intervals $[u_k - \alpha_k^{-1}, u_k + \alpha_k^{-1}]$ overlap.

The signal that our wavelet estimator will endeavor to recover is

(4.10) $$g = g_0 + \sum_k \gamma\{\alpha_k(\cdot - u_k)\}.$$

For the present we assume that the interval $\mathcal{J}$ is placed immediately to the right of $\mathcal{I}$, so that (in view of Theorem 4.1) the probability that the algorithm enters $\mathcal{J}$ using sampling rate $\rho_1$ converges to 1 as $\nu \to \infty$. Our next result gives conditions under which, if the fluctuations in $\mathcal{J}$ are as described at (4.10), then (with high probability) a rate switch from $\rho_1$ to $\rho_2$ occurs during $\mathcal{J}$.

If the frequency $\alpha_1$ of the first fluctuation satisfies (4.6) for $k = 1$, then, with probability converging to 1 as $\nu \to \infty$, there will be a switch to rate $\rho_2$ in the close vicinity of time $u_1$. This follows from Theorem 4.4, on considering the case where the series at (4.10) consists of the single fluctuation $\gamma\{\alpha_1(\cdot - u_1)\}$. In such a case, the theorem does not make any comment on what happens later in interval $\mathcal{J}$; that will be dealt with in Theorems 4.5 and 4.6.

THEOREM 4.4. *Suppose data are generated by the model at (2.1), with signal $g$ given by (4.10). Assume the estimator $\hat{g}$ is constructed using $C > 2^{1/2}$ in the threshold $\delta$, that the rule* (R.1) *is used to define an upward rate switch, and that (4.1)–(4.4) hold, with $B > C^2/(1-c)$ in (4.2) and, in (4.3), $(q, \rho) = (q_1, \rho_1)$. Suppose too that (4.6)–(4.9) hold. If, on entering time interval $\mathcal{J}$, the sampling rate is $\rho_1$, then with probability converging to 1 as $\nu \to \infty$, an increase in the rate to rate $\rho_2$ will occur during time interval $\mathcal{J}$.*

We continue to assume the signal is composed of fluctuations that may be modelled as at (4.10). However, when showing that the rate will not change during the time interval $\mathcal{J}$, we make the additional assumption that during each subinterval of $\mathcal{J}$, of length $\tau_0$, there exists a fluctuation whose frequency is of larger order than $\pi_2$ and of smaller order than $p_{q_2}$:

> it is possible to choose a subset $\mathcal{A}$ of the set of all frequencies $\alpha_k$ represented at (4.10), such that, for each interval $\mathcal{K}$ of
(4.11) length $\tau_0$ included within $\mathcal{J}$, there is at least one $\alpha_k \in \mathcal{A}$ such that the associated function $\gamma\{\alpha_k(\cdot - u_k)\}$ is supported within $\mathcal{K}$, and, moreover, $\pi_2 = o(\min_{\alpha \in \mathcal{A}} \alpha)$ and $\max_{\alpha \in \mathcal{A}} \alpha = o(p_{q_2})$.

In the result below, we assume that we start the time interval $\mathcal{J}$ using sampling at rate $\rho_2$. Thus, $\mathcal{J}$ can no longer be thought of as following immediately after an interval where the signal is smooth. However, it could follow immediately after a short interval that contained a single fluctuation $\alpha = \alpha_1$ which triggered a switch from rate $\rho_1$ to $\rho_2$; see the paragraph immediately preceding Theorem 4.4. Recall that rule (R.2), for switching to a lower sampling rate, was given in Section 2.

THEOREM 4.5. *Assume that the estimator $\hat{g}$ is constructed using $C > 2^{1/2}$ in the threshold $\delta$, that the rule* (R.2) *is used to define a downward rate switch, and that (4.1)–(4.4) hold, with $B > C^2/(1-c)$ in (4.2) and, in (4.3),* $(q, \rho) = (q_2, \rho_2)$. *Suppose too that* $|\mathcal{J}|$ *is bounded as* $\nu \to \infty$, *that* (4.7)–(4.9) *and* (4.11) *hold, that* $g$ *is given by* (4.10), *and that the sampling rate at the start of* $\mathcal{J}$ *equals* $\rho_2$. *Then, with probability converging to* 1 *as* $\nu \to \infty$, *the sampling rate stays at* $\rho_2$ *throughout* $\mathcal{J}$.

This result has an analogue in which the frequencies in $\mathcal{J}$ are relatively low, and a switch from sampling rate $\rho_2$ to $\rho_1$ is virtually assured:

THEOREM 4.6. *Assume the conditions in Theorem* 4.5, *except that the constraints* "$\pi_2 = o(\min_{\alpha \in \mathcal{A}} \alpha)$ *and* $\max_{\alpha \in \mathcal{A}} \alpha = o(p_{q_2})$" *at the end of* (4.11) *are changed to* "$\max_{\alpha \in \mathcal{A}} \alpha = o\{\min(\pi_1, \pi_2)\}$." *Then, with probability converging to* 1 *as* $\nu \to \infty$, *the sampling rate switches from* $\rho_2$ *to* $\rho_1$ *during* $\mathcal{J}$, *and stays there for the duration of that time interval.*

Finally we show that, when sampling is carried out at rate $\rho$, the estimator is able to consistently recover frequencies almost up to the level $\rho$.

THEOREM 4.7. *Suppose data are generated by the model at* (2.1), *with independent errors* $\varepsilon_i$. *Assume the sampling rate is constant at* $\rho = \rho_1$ *or* $\rho_2$, *and that the threshold* $\delta$ *is given by* (2.4) *with* $C > 2^{1/2}$. *Suppose too that* (4.2)–(4.4) *hold, with* $B > C^2/(1 - c)$ *in* (4.2) *and, in* (4.3), $(q, \rho) = (q_1, \rho_1)$ *or* $(q_2, \rho_2)$ *for the respective sampling rates. Assume the signal is given by* (4.10) *on* $\mathcal{J}$, *where* $\max \alpha_k = o(p_q)$. *Then, for each* $\eta > 0$, *the probability that* $|\hat{g} - g| \leq \eta$ *uniformly on* $\mathcal{J}$ *converges to* 1 *as* $\nu \to \infty$.

We conclude by quantifying some of the potential gains and losses from dual-rate sampling. Suppose the expense of sampling, expressed, for example, in terms of the capacity of the data storage device, demands that the long-run sampling rate not exceed $\nu$ per unit time. If, in parts of the signal that have relatively high frequency, we use rate $\rho_2 > \nu$ rather than $\nu$, then (in order to stay within budget) at other time points we should reduce the rate to $\rho_1$, where $\rho_1$ and $\rho_2$ are connected by the formula

$$(4.12) \qquad\qquad \nu = \rho_1(1 - \Pi) + \rho_2 \Pi,$$

and $\Pi$ denotes the long-run proportion of time for which we use rate $\rho_2$.

It may be deduced from Theorem 4.2 that the condition for there to be no asymptotic deterioration in mean-squared error, to first order, in the relatively smooth places where rate $\rho_1$ is employed, is $\nu \sim \rho_1$. This is, of course, equivalent to $\Pi \rho_2 \to 0$ as $\nu \to \infty$. In the proportion $\Pi$ of the time when we use the higher sampling rate, there is (in view of Theorem 4.7) potential for consistently estimating the signal where this would not otherwise be possible.

# REFERENCES

ALDROUBI, A. and GRÖCHENIG, K. (2001). Nonuniform sampling and reconstruction in shift-invariant spaces. *SIAM Rev.* **43** 585–620. MR1882684

ANON (2002). Mathematics: A new wave. *The Economist* January 19 70. MR1913633

ANTONIADIS, A. and FAN, J. Q. (2001). Regularization of wavelet approximations (with discussion). *J. Amer. Statist. Assoc.* **96** 939–967. MR1946364

CAI, T. T. and BROWN, L. D. (1998). Wavelet shrinkage for nonequispaced samples. *Ann. Statist.* **26** 1783–1799. MR1673278

CAI, T. T. and BROWN, L. D. (1999). Wavelet estimation for samples with random uniform design. *Statist. Probab. Lett.* **42** 313–321. MR1688134

CHEN, W., ITOH, S. and SHIKI, J. (1998). Irregular sampling theorems for wavelet subspaces. *IEEE Trans. Inform. Theory* **44** 1131–1142. MR1616719

DELOUILLE, V., FRANKE, J. and VON SACHS, R. (2001). Nonparametric stochastic regression with design-adapted wavelets. *Sankhyā Ser. A* **63** 328–366. MR1897046

DELOUILLE, V., SIMOENS, J. and VON SACHS, R. (2001). Smooth design-adapted wavelets for nonparametric stochastic regression. Discussion Paper No. 0117, l'Institut de Statistique, Univ. Catholique de Louvain, Belgium.

DELOUILLE, V. and VON SACHS, R. (2002). Properties of design-adapted wavelet transforms of nonlinear autoregression models. Discussion Paper No. 0225, l'Institut de Statistique, Univ. Catholique de Louvain, Belgium.

DONOHO, D. L. and JOHNSTONE, I. M. (1994). Ideal spatial adaptation by wavelet shrinkage. *Biometrika* **81** 425–455. MR1311089

DONOHO, D. L. and JOHNSTONE, I. M. (1995). Adapting to unknown smoothness via wavelet shrinkage. *J. Amer. Statist. Assoc.* **90** 1200–1224. MR1379464

DONOHO, D. L., JOHNSTONE, I. M., KERKYACHARIAN, G. and PICARD, D. (1995). Wavelet shrinkage: Asymptopia? (With discussion.) *J. Roy. Statist. Soc. Ser. B* **57** 301–369. MR1323344

HALL, P. and PENEV, S. (2002). Wavelet-based estimation with multiple sampling rates. Available at www.maths.unsw.edu.au/~spiro/publicat.html.

HALL, P. and PATIL, P. (1995). Formulae for mean integrated squared error of nonlinear wavelet-based density estimators. *Ann. Statist.* **23** 905–928. MR1345206

HALL, P., PARK, B. U. and TURLACH, B. A. (1998). A note on design transformation and binning in nonparametric curve estimation. *Biometrika* **85** 469–476. MR1649126

HALL, P. and TURLACH, B. A. (1997). Interpolation methods for nonlinear wavelet regression with irregularly spaced design. *Ann. Statist.* **25** 1912–1925. MR1474074

HERZBERG, A. M. and TRAVES, W. N. (1994). An optimal experimental design for the Haar regression model. *Canad. J. Statist.* **22** 357–364. MR1309320

KOVAC, A. and SILVERMAN, B. W. (2000). Extending the scope of wavelet regression methods by coefficient-dependent thresholding. *J. Amer. Statist. Assoc.* **95** 172–183.

OYET, A. J. and WIENS, D. P. (2000). Robust designs for wavelet approximations of regression models. *J. Nonparametr. Statist.* **12** 837–859. MR1802579

PENSKY, M. and VIDAKOVIC, B. (2001). On non-equally spaced wavelet regression. *Ann. Inst. Statist. Math.* **53** 681–690. MR1879604

SARDY, S., PERCIVAL, D. B., BRUCE, A. G., GAO, H. Y. and STUETZLE, W. (1999). Wavelet shrinkage for unequally spaced data. *Statist. Comput.* **9** 65–75.

ZHANG, S. L. and ZHENG, Z. G. (1999). Nonlinear wavelet estimation of regression function with random design. *Sci. China Ser. A* **42** 825–833. MR1738553

CENTRE FOR MATHEMATICS
  AND ITS APPLICATIONS
AUSTRALIAN NATIONAL UNIVERSITY
CANBERRA, ACT 0200
AUSTRALIA
AND
DEPARTMENT OF STATISTICS
LONDON SCHOOL OF ECONOMICS
HOUGHTON STREET
LONDON WC2A 2AE
UNITED KINGDOM
E-MAIL: Peter.Hall@maths.anu.edu.au

DEPARTMENT OF STATISTICS
SCHOOL OF MATHEMATICS
UNIVERSITY OF NEW SOUTH WALES
2052 SYDNEY, NSW
AUSTRALIA
E-MAIL: spiro@maths.unsw.edu.au